\documentclass[12pt]{amsart}

\usepackage{amssymb, amsmath,amsthm}
\usepackage{slashbox}
\usepackage{MnSymbol}
\usepackage[english]{babel}
\usepackage[utf8]{inputenc}
\newtheorem{thm}{Theorem}
\newtheorem{lem}[thm]{Lemma}
\newtheorem{cor}[thm]{Corollary}
\newtheorem{rem}{Remark}

\begin{document}
\title[The $2$-adic valuation of generalized Fibonacci sequences]{The 2-adic valuation of generalized Fibonacci sequences with an application to certain Diophantine equations.}
\author{Bartosz Sobolewski}
\address{Jagiellonian University, Faculty of Mathematics and Computer Science, Institute of Mathematics, {\L}ojasiewicza 6, 30 - 348 Krak\'{o}w, Poland}
\email{\tt bartosz.sobolewski@doctoral.uj.edu.pl}

\begin{abstract}
In this paper we focus on finding all the factorials expressible as a product of a fixed number of $2k$-nacci numbers with $k \geq 2$. We derive the 2-adic valuation of the $2k$-nacci sequence and use it to establish bounds on the solutions of the initial equation.
In addition, we specify a more general family of sequences, for which we can perform a similar procedure. We also investigate a possible connection of these results with $p$-regular sequences.
\end{abstract}
\maketitle
\section{Introduction} \label{sec:introduction}

For a fixed integer $r \geq 2$  define the generalized Fibonacci ($r$-nacci) sequence $\{t_n\}_{n\geq 0}$ as follows:
\begin{align}\label{def:2k-nacci}
t_n =
\begin{cases}
0 & \text{ for } n = 0,\\ 
1 & \text{ for } 1 \leq n \leq r-1, \\ 
\sum\limits_{i=1}^{r} t_{n-i} & \text{ for } n \geq r-1.
\end{cases}
\end{align}
Notice, that for $r=2$ we obtain the usual Fibonacci sequence, which has already been studied extensively. In this paper we will mostly restrict ourselves to the case of even $r \geq 4$ and write $r=2k$ for some $k \geq 2$.
The main motivation for our considerations is to completely solve the equation
\begin{align} \label{eq:product_of_terms}
m! = \prod_{i=1}^d t_{n_i}
\end{align}
in positive integers $m, n_1, ..., n_d$.

For $p$ prime define the $p$-adic valuation of a non-zero integer $s$ as $\nu_p(s) = \max \{l \geq 0 : p^l | s \}$ and $\nu_p(0) = \infty$. 
Equation (\ref{eq:product_of_terms}) for the case of $r=3$ and $d=1$ was solved by Lengyel and Marques in \cite{Lengyel_Marques} by means of computing $\nu_2(t_n)$ and then applying this result to obtain an effective upper bound on $m$ and $n_1$. In this paper we will follow a similar procedure for $2k$-nacci sequences with $k \geq 2$.

To begin with, in Theorem \ref{thm:equation_solution} we specify a more general family of integer sequences $\{s_n\}_{n\geq 0}$ for which we are able to solve equation (\ref{eq:product_of_terms}) and show a general procedure to achieve this goal. Informally speaking, we need the term $s_n$ to grow at least exponentially and $\nu_p(s_n)$ -- at most polynomially with an exponent less than 1, for some $p$ prime.

Theorem \ref{thm:2-adic_order} provides a simple expression for $\nu_2(t_n)$ when $r=2k \geq 4$ and the subsequent corollary shows that the sequence $\{t_n\}_{n\geq 0}$ satisfies the conditions of Theorem \ref{thm:equation_solution}. We then find all the solutions of equation (\ref{eq:product_of_terms}) for small values of $k$ and $d$.

We also briefly discuss how our results are related to $p$-regular sequences. Recall that a sequence $\{s_n\}_{n\geq 0}$ with rational values is $p$-regular iff its $p$-kernel
$$
\mathcal{N}_p(a) = \left\{ \{s_{p^l n + j} \}_{n\geq 0} : \; l\geq 0, \; 0 \leq j < p^l\right\}
$$
is contained in a finitely generated $\mathbb{Z}$-module. More details on regular sequences can be found in \cite{Allouche_Shallit} and \cite{Allouche_Shallit_2}.
As we note later, the formula given in Theorem \ref{thm:2-adic_order} implies $2$-regularity of $\{\nu_2(t_n)\}_{n \geq 0}$. However, it turns out that exponential growth of a sequence $\{s_n\}_{n \geq 0}$ and $p$-regularity of $\{\nu_p(s_n)\}_{n \geq 0}$ still do not guarantee that the assumptions of Theorem \ref{thm:equation_solution} are met. In this case we cannot determine, using the shown method, whether the equation (\ref{eq:product_of_terms}) has only a finite number of solutions.

\section{Main results}\label{sec:main}
As we mentioned before, we start with describing a general situation in which equation (\ref{eq:product_of_terms}) can be completely solved.
For two sequences $\{a_n\}_{n \geq 0}$ and $\{b_n\}_{n \geq 0}$ we denote $a_n = O(b_n)$ if there exists a positive constant $K$ such that $|a_n| \leq K |b_n|$ for sufficiently large $n$. Similarily, we write $a_n = \Omega(b_n)$ if there exists a positive constant $K$ such that $|a_n| \geq K |b_n|$ for sufficiently large $n$.
First, we give an auxiliary lemma, also used in \cite{Lengyel_Marques}, which is an easy corollary from Legendre's formula for $\nu_p(m!)$. 

\begin{lem}\label{lem:legendre_ineq}
For any integer $m \geq 1$ and prime $p$, we have
$$\frac{m}{p-1} - \left \lfloor \frac{\log m}{\log p} \right \rfloor -1 \leq \nu_p(m!) \leq \frac{m-1}{p-1}.$$
\end{lem}

\begin{thm} \label{thm:equation_solution}
Let $\{s_n\}_{n\geq 0}$ be a sequence of positive integers such that 
\begin{align} \label{eq:from_below}
\log s_n = \Omega(n).
\end{align} 
Let $p$ be a prime. Assume that
\begin{align} \label{eq:from_above}
\nu_p(s_n) = O\left(n^C\right)
\end{align}
for some constant $C < 1$. Then for each fixed positive integer $d$ the equation
\begin{align} \label{eq:general_product}
m! = \prod_{i=1}^d s_{n_i}
\end{align}
has only a finite number of solutions in $m, n_1, n_2, \ldots, n_d$.
\begin{proof}
We adjust the method used in \cite{Lengyel_Marques} to a more general setting. Roughly speaking, we will show that if
 (\ref{eq:general_product}) is satisfied and we let both sides grow, then the $p$-adic valuation of the right hand side increases slower than $\nu_p(m!)$. For each value of $p$ we proceed in the same way, so for simplicity assume that $p=2$.
By our assumptions, there exist some positive constants $K_1, K_2$ and an integer $n_0 \geq 0$ such that $\nu_2(s_n) \leq K_1 n^C$ and $\log_2 s_n \geq K_2 n$ for $n \geq n_0$. 

Suppose that $n_i \geq n_0$ for $i=1,2,\ldots,d$. There is only a finite number of solutions with $m < 6$ because $s_n$ grows at least exponentially.
By Lemma \ref{lem:legendre_ineq}, for $p = 2$ we get 
\begin{align} \label{eq:m_lower_bound}
\frac{1}{2}m \leq m - \lfloor \log_2 m \rfloor - 1  \leq \nu_2(m!), 
\end{align}
where the leftmost inequality is true for $m \geq 6$. On the other hand,
\begin{align*}
\nu_2 \left( \prod_{i=1}^d s_{n_i} \right ) = \sum_{i=1}^{d} \nu_2(s_{n_i}) \leq K_1 \sum_{i=1}^{d} n_i^C
\leq d K_1 \left(\max_{1 \leq i \leq d} n_i \right)^C.
\end{align*}
Hence, for $m \geq 6$
\begin{align} \label{eq:m_below_n}
m \leq 2dK_1 \left(\max_{1 \leq i \leq d} n_i \right)^C.
\end{align}
We need another inequality with $n_i$ bounded from above by $m$. By the assumption (\ref{eq:from_below})
\begin{align*}
\log_2 \left(\prod_{i=1}^d s_{n_i} \right) = \sum_{i=1}^{d} \log_2 s_{n_i} \geq K_2 \sum_{i=1}^{d} n_i \geq  K_2 \max_{1 \leq i \leq d} n_i.
\end{align*}
Furthermore, for $m \geq 5$ the following inequality holds:
\begin{align} \label{eq:m_fact_from_above}
m! < \left ( \frac{m}{2} \right)^m,
\end{align}
and hence
\begin{align} \label{eq:n_below_m}
 m \log_2 \frac{m}{2} > K_2 \max_{1 \leq i \leq d} n_i.
\end{align}
Combining inequalities (\ref{eq:m_below_n}) and (\ref{eq:n_below_m}) yields
\begin{align} \label{eq:m_below_m}
\log_2 \frac{m}{2} > K m^{\frac{1}{C} -1},
\end{align}
where $K = K_2/ (2dK_1)^{1/C}$. But $C < 1$ implies that the exponent in (\ref{eq:m_below_m}) is strictly greater than 0, so (\ref{eq:m_below_m}) holds only for finitely many $m$. Thus, by inequality (\ref{eq:n_below_m}) we obtain an upper bound for all $n_i$.

Now assume without loss of generality that $n_i < n_0$ for $i=j+1,j+2,\ldots,d$. Observe that for each fixed $n_{j+1},n_{j+2},\ldots,n_d$ the problem is equivalent to solving the equation (\ref{eq:general_product}) with $d=j$ and the left hand side divided by a positive integer constant $M$.
Then in (\ref{eq:m_lower_bound}) and (\ref{eq:m_fact_from_above}) we need to replace $m!$ with $m! /M$ which only changes the set of $m$ for which both of those inequalities hold.
But this leads to the same conclusion as before.
\end{proof}
\end{thm}

\begin{rem}
Observe that the method of Theorem \ref{thm:equation_solution} does not work if we let $d$ be unbounded. Indeed, the constant $K$ in (\ref{eq:m_below_m}) becomes arbitrarily small as $d$ increases, so we cannot use the subsequent argument. 
Informally speaking, the $p$-adic valuation of the expression $\prod_{i=1}^d s_{n_i}$ might grow too fast for the method to work.
\end{rem}

\begin{rem}
The condition (\ref{eq:from_below}) is satisfied for sequences expressible in Binet form. Hence, for linear recurrence sequences, we usually need to check only the condition (\ref{eq:from_above}) for some $p$. 
One might also ask whether we can replace it with some other assumption.
Shu and Yao proved in \cite{Shu_Yao} a condition on a binary recurrence sequence, which guarantees $p$-regularity of $\{\nu_p(s_n)\}_{n \geq 0}$, and mentioned a possible generalization to recurrences of higher order.  It is known that $p$-regular sequences grow at most polynomially, which is a result by Allouche and Shallit \cite{Allouche_Shallit}.
Unfortunately, this does not give a bound on $C$ in (\ref{eq:from_above}) and the proof of Theorem \ref{thm:equation_solution} fails if $C \geq 1$. 
Therefore, some additional information besides regularity needs to be known about $\{\nu_p(s_n)\}_{n \geq 0}$ in order to put the theorem to use.
\end{rem} 

The reasoning in Theorem \ref{thm:equation_solution} provides an upper bound on the solutions of the equation (\ref{eq:general_product}) if we are able to find the values of $C,K_1, K_2$ and $n_0$. We will show that it is indeed the case for our sequence $\{t_n\}_{n \geq 0}$. 
We start by determining the $2$-adic valuation of each $t_n$. A similar characterization of $\{\nu_2(t_n)\}_{n\geq 0}$ for $r=2$ is given by Lengyel in \cite{Lengyel_order}  and \cite{Lengyel_multi} (by a different method).

\begin{thm}\label{thm:2-adic_order}
If $k \geq 2$, then the sequence $\{\nu_2(t_n)\}_{n\geq 0}$ satisfies the following conditions:
$$
\nu_2(t_n) = 
\begin{cases}
0 & \text{ for }	n \equiv 1,2,...,2k \pmod{2k+1}, \\
1 & \text{ for }	n \equiv 2k+1		\pmod{2(2k+1)}, \\
\nu_2(n) + \nu_2(k-1) + 2 & \text{ for }	n \equiv 0 \pmod{2(2k+1)}.	
\end{cases}
$$
\end{thm}

The proof of the theorem is presented in Section \ref{sec:proof}. 

\begin{rem}
Theorem \ref{thm:2-adic_order} implies that $\{\nu_2(t_n)\}_{n\geq 0}$ is a $2$-regular sequence. Note that this conclusion does not follow directly from the results of \cite{Shu_Yao}, as we consider recurrence of any even order.
\end{rem}

Now we proceed to show that Theorem \ref{thm:equation_solution} can be applied to our sequence $\{t_n\}_{n \geq 0}$.
The following lemma establishes a lower bound on $t_n$.

\begin{lem} \label{lem:term_size}
For all $n \geq 1$ we have
\begin{align} \label{eq:exp_growth}
t_n \geq \phi^{n-r-1},
\end{align}
where $\phi$ is the unique real root of the equation $x^{r} = x^{r-1} + ... + x + 1$ lying in the interval $(1,2)$.
\begin{proof}
By lemma 3.6 in \cite{Wolfram} there is indeed exactly one such $\phi$, which in addition lies in the interval $(2(1- 2^{-r}), 2)$. For $n = 1, 2, ..., r$ the inequality (\ref{eq:exp_growth}) follows from starting conditions for $t_n$ and the fact that $t_r = 2(r-1) \geq 2 > \phi$. Then we proceed easily by induction.
\end{proof}
\end{lem}

\begin{cor}\label{cor: finite_solutions}
If $r=2k \geq 4$ then the equation (\ref{eq:product_of_terms}) has only a finite number of solutions in $m, n_1, n_2, \ldots, n_d$ and this number can be effectively bounded from above.
\begin{proof}
By theorem \ref{thm:2-adic_order} we have
$$\nu_2(t_n) \leq \nu_2(n) + \nu_2(k-1) + 2 \leq \log_2(n) + \nu_2(k-1) + 2 \leq \sqrt{n},$$
where the last inequality is true for example for $n \geq 2^{2 \max \{2, \nu_2(k-1)\}}$.
By lemma \ref{lem:term_size}
$$ \log t_n \geq (n-2k-1) \log \phi \geq \frac{1}{2} n \log \phi,$$
for $n \geq 2(2k+1)$, because $\phi > 1$. Hence, we can take $n_0 = \max \{2(2k+1), 2^{2 \max \{2, \nu_2(k-1)\}}\}$ and apply the method used in Theorem \ref{thm:equation_solution} to find an upper bound on $m, n_1, n_2, \ldots, n_d$.
\end{proof}
\end{cor}

We could apply Corollary \ref{cor: finite_solutions} to find all non-trivial solutions (with $t_{n_i} > 1$ for each $i$)  of the equation (\ref{eq:product_of_terms}) for given $k$ and $d$.
 However, using the explicit form of $\nu_2(t_n)$, one can make the bounds much more precise. We will once again follow the approach shown in \cite{Lengyel_Marques}. The computations are quite similar to those in Theorem \ref{thm:equation_solution}, so we omit the details.

By  Lemma \ref{lem:legendre_ineq} and Theorem \ref{thm:2-adic_order}, we obtain
\begin{align}\label{eq:m_below_n_2}
m - \lfloor \log_2 m \rfloor - 1 &\leq \nu_2(m!) =  \sum_{i=1}^{d} \nu_2(t_{n_i}) \leq \sum_{i=1}^{d} [\nu_2(n_i) + \nu_2(k-1) + 2 ] \nonumber \\ &= d[\nu_2(k-1) + 2] + \nu_2 \left( \prod_{i=1}^d n_i \right ) \nonumber \\
&\leq d[\nu_2(k-1) + 2] + \log_2 \left( \prod_{i=1}^d n_i \right ) .
\end{align}
On the other hand, from Lemma \ref{lem:term_size}, and inequality (\ref{eq:m_fact_from_above}),  we get for $m \geq 5$
\begin{align}\label{eq:n_below_m_2}
\left [\sum_{i=1}^{d} n_i - d(2k-1) \right ] {\log_2 \phi_k} \leq  \sum_{i=1}^{d} \log_2(t_{n_i}) = \log_2(m!) < m (\log_2 m -1),
\end{align}
where $\phi_k$ is the value of $\phi$ in Lemma \ref{lem:term_size} corresponding to $r=2k$.
The AM--GM inequality applied to all $n_i$, together with  (\ref{eq:m_below_n_2}) and (\ref{eq:n_below_m_2}), yields
\begin{align}\label{eq:only_m}
m - \lfloor \log_2 m \rfloor - 1 - d[\nu_2(k-1) + 2] < d \log_2
\left [ \frac{m}{d\log_2 \phi_k} (\log_2 m -1) + (2k-1)   \right ].
\end{align}
This gives an upper bound on $m$ and, consequently, on each $n_i$. As an example, in the table below we give the upper bound on $m$ obtained for $2k$-nacci sequences with $2 \leq k \leq 5$ and $1 \leq d \leq 10$.
\begin{table}[!htbp]
\centering
\label{my-label}
\begin{tabular}{c|cccccccccc}
\backslashbox{$k$\kern-0.5em}{\kern-0.5em $d$}   & 1  & 2  & 3  & 4  & 5  & 6  & 7  & 8  & 9  & 10   \\ \hline
2 & 11 & 19 & 27 & 35 & 43 & 51 & 59 & 67 & 75 & 84  \\
3 & 13 & 22 & 31 & 40 & 50 & 59 & 68 & 77 & 87 & 96  \\
4 & 11 & 19 & 28 & 36 & 44 & 52 & 60 & 68 & 76 & 84  \\
5 & 14 & 25 & 35 & 46 & 56 & 67 & 77 & 88 & 98 & 109
\end{tabular}
\end{table}

Using this result we find that the only non--trivial solution of the equation (\ref{eq:product_of_terms}) with $2 \leq k \leq 5$ and $1 \leq d \leq 10$ appears in the $4$-nacci sequence and is the single term $t_5 = 3!$.

\section{Proof of Theorem \ref{thm:2-adic_order}} \label{sec:proof}
In order to study the 2-adic valuation of $t_n$ it is enough to focus on $n$ divisible by $2k+1$ as all the other terms are odd. In this case we can write $n = s 2^l (2k+1)$ for $s$ odd and $l \geq 0$.
We will divide our proof into two main parts.

First, we will show by induction on $l$ and $s$ that $2k$ consecutive terms of the sequence $\{t_n\}_{n \geq 0}$, starting with $t_{s 2^l (2k+1)}$, satisfy a particular system of congruences, given in Lemma \ref{lem:later_terms}.
However, as it will turn out, this argument works only for $l \geq l_0$, where $l_0$ depends on $\nu_2(k-1)$.
Moreover, the initial system of congruences for $l=l_0$ involves some constants, which need to be computed.
 
 Therefore, we will have to employ another method for $l \leq l_0$. We will show how to obtain the values of $t_{n + 2k+1}, \ldots, t_{n+ 4k}$ in terms of  $t_{n}, \ldots, t_{n+ 2k-1}$ for any $n$ and then proceed by induction. 
 As a result, we will be able to express $t_n$ in quite concrete form, given in Lemma \ref{lem:easy_expression}, involving binomial coefficients weighted by powers of 2.

For simplicity, we introduce the following matrix notation:
\begin{align*}
T_n = \begin{bmatrix}
t_n \\
t_{n+1} \\
\vdots \\
t_{n+2k-1}
\end{bmatrix},
\qquad
B_n =
\begin{bmatrix}
t_n & t_{n+1} & \hdots & t_{n+ 2k-1} \\ 
t_{n+1} & t_{n+2} & \hdots & t_{n+2k}\\ 
\vdots& \vdots &  &\vdots\\ 
t_{n+ 2k-1} &t_{n+2k} & \hdots & t_{n+4k-2}
\end{bmatrix},
\end{align*}
where $n \geq 0$. By $C$ we will denote the companion matrix of $t_n$ which has the form
\begin{align} \label{eq:matrix_C}
C=
\begin{bmatrix}
0 & 1 & 0 &\hdots & 0 & 0 \\ 
\vdots & 0 & 1 & \ddots &  & 0\\ 
 & \vdots & \ddots &\ddots & \ddots & \vdots \\ 
 &  &  & \ddots & 1 & 0	\\
0 & 0 & 0 & \cdots & 0 & 1 \\
1 & 1 & 1 &\hdots & 1 & 1
\end{bmatrix},
\end{align}
where the entries above the diagonal and in the bottom row are equal to 1 and all other entries are zero.
It is easy to check that $CT_n = T_{n+1}$ and $CB_n = B_{n+1}$, so for any positive integers $n$ and $w$ we have
\begin{align}
C^n T_w &= T_{n+w}, \label{eq:pushing_the_vector} \\
C^n B_w &= B_{n+w}. \label{eq:pushing_the_matrix}
\end{align}

First, we state an identity involving the terms of the sequence $\{t_n\}_{n \geq 0}$.

\begin{lem} \label{lem:reduction_formula}
The matrix $B_0$ is invertible and its determinant is odd. Moreover,
for all positive integers $n, w$ we have the formula
\begin{align} \label{eq:reduction_formula}
t_{n+w} = T_n^T B_0^{-1} T_w.
\end{align}

\begin{proof}
It is easy to see that $t_n$ is even iff $n$ is divisible by $2k+1$. Therefore,
$$
\det B_0 \equiv \det
\begin{bmatrix}
0 & 1 &1 & \hdots & & 1 \\
1 & 1 & \hdots  &&& 1 \\
1 & \vdots& &&\udots & 0 \\
\vdots&& & \udots & \udots & 1 \\
 & &\udots & \udots & \udots & \vdots \\
1 & 1 &  0 & 1 & \cdots & 1
\end{bmatrix} \pmod{2},
$$
where zeros in the latter matrix appear only at positions corresponding to $t_0$ and $t_{2k+1}$ in $B_0$, that is, at $(1,1)$ and $(i,j)$ such that $i+j = 2k + 2$. By subtracting the first row from all the others, we easily obtain $\det B_0 \equiv 1 \pmod{2}$, which proves the first part of the statement.

From (\ref{eq:pushing_the_vector}) and (\ref{eq:pushing_the_matrix}), we get
$$
T_{n+w} = C^n T_{w} = B_n B_0^{-1} T_w.
$$
The first coordinate gives us the formula for $t_{n+w}$.
\end{proof}
\end{lem}

The identity (\ref{eq:reduction_formula}) might seem difficult to apply without an explicit expression for $B_0^{-1}$. However, it plays a major role in deriving the congruence relations in the following lemma.

\begin{lem} \label{lem:later_terms}
For any $l \geq 0$ the following congruence relation holds: 
\begin{align} \label{eq:easy_divis}
T_{2^l (2k+1)} \equiv T_0 \pmod{2^{l+1}}.
\end{align}
Moreover, if a column vector $A \in \mathbb{Z}^{2k}$ satisfies
\begin{align*} 
T_{2^{l_0} (2k+1)} \equiv 2^{l_0+1} A + T_0  \pmod{2^{l_0+\nu_2(k-1)+3}},
\end{align*}
where $l_0 = \nu_2(k-1)+2$, then for any $l \geq l_0$ and $s \geq 1$ we also have
\begin{align} \label{eq:the_second_form_2}
T_{s2^l (2k+1)} \equiv s2^{l+1} A + T_0  \pmod{2^{l+\nu_2(k-1)+3}}.
\end{align}
\begin{proof}
Obviously (\ref{eq:easy_divis}) is true for $l=0$. 
Now assume that  (\ref{eq:easy_divis}) holds for some $l \geq 0$. We can write
\begin{align} \label{eq:the_second_form}
t_{2^l (2k+1)+ j} = 2^{l+1} a_{l,j} + t_j,
\end{align}
where $a_{l,j}$ are some positive integers for $j=0,1, \ldots, 2k-1$. Define also $a_{l,2k}, a_{l,2k+1}, \ldots,a_{l,4k-2}$ by the same recurrence  as $\{t_n\}_{n \geq 0}$. Then (\ref{eq:the_second_form}) is satisfied for $j=0,1, \ldots, 4k-2$.  For convenience denote by $e_j \in \mathbb{Z}^{2k}$ the vector with 1 on the $j$-th position (counting from 0) and 0 on the other positions, and additionally define
\begin{align*} 
A_{l,j} = 
\begin{bmatrix}
a_{l,j} &
a_{l,j+1} &
\cdots &
a_{l,j+2k-1}
\end{bmatrix}^T
\end{align*}
for $j=0,1, \ldots, 2k-1$. It follows from the definition of $B_0$  that $B_0^{-1} T_j = e_j$.

Fix any $0 \leq i \leq 2k-1$. The formula  (\ref{eq:reduction_formula}) yields
\begin{align*} 
t_{2^{l+1}(2k+1) + i} = T_{2^{l}(2k+1)}^T B_0^{-1} T_{2^{l}(2k+1) + i}
\end{align*}
Therefore, using (\ref{eq:the_second_form}) we get
\begin{align} \label{eq:reduction_6}
t_{2^{l+1}(2k+1) + i} &= 2^{2l+2} A_{l,0}^T B_0^{-1} A_{l,i} + 2^{l+1} \left( A_{l,0}^T e_i +e_0^T A_{l,i}
\right ) + T_0^T e_i \nonumber \\
 &= 2^{2l+2}c_{l,i} + 2^{l+2} a_{l,i} + t_i
\end{align}
for some rational $c_{l,i}$ such that $\det(B_0) c_{l,i}$ is an integer. In Lemma \ref{lem:reduction_formula}, however, we showed that $\det B_0$ is odd which means that $c_{l,i}$ must be an integer. Thus,
\begin{align*} 
2^{l+2} |t_{2^{l+1}(2k+1) + i} - t_i,
\end{align*}
from which (\ref{eq:easy_divis}) follows. If we choose $l \geq l_0 = \nu_2(k-1)+2$ then the term $2^{2l+2}  c_{l,i}$ in (\ref{eq:reduction_6}) is reduced modulo $2^{l+\nu_2(k-1)+4}$. We can take $A = A_{l_0,0}$ to complete the proof of (\ref{eq:the_second_form_2}) for $s=1$. 

To proceed by induction on $s$ notice that the index $s 2^l (2k+1) + i$ can be expressed as a sum of indices in the following way:
$$s 2^l (2k+1) + i = (2^l (2k+1) + i) + (s-1) 2^l (2k+1).$$ 
We can then perform a similar computation as in (\ref{eq:reduction_6}) to get the desired result.
\end{proof}
\end{lem}

Our specific choice of the divisor in (\ref{eq:the_second_form_2}) equal to $2^{l+\nu_2(k-1)+3}$ is based on the observation of $t_n$ and is indeed effective in proving the formula for $\nu_2(t_n)$. The numbers $a_{l,j}$ in (\ref{eq:the_second_form}) are determined uniquely modulo $2^{\nu_2(k-1)+2}$. We are particularly interested in finding the value of $a_{l_0,0}$ which will directly give us the 2-adic valuation of $t_{2^l(2k+1)}$ for $l \geq l_0 =  \nu_2(k-1)+2$, provided that $a_{l_0,0} \leq \nu_2(k-1)+1$. However, as mentioned before, we need to develop another method to analyze the case when $l \leq l_0$. 

We start with deriving a formula for expressing $t_{n + 2k+1}, \ldots, t_{n+ 4k}$ in terms of  $t_{n}, \ldots, t_{n+ 2k-1}$.

\begin{lem}\label{lem:pushing_the_sequence}
Define $C$ as in (\ref{eq:matrix_C}). Then
$$
C^{2k+1} = 2
\begin{bmatrix}
1 & 1 & \hdots & 1 \\
2 & 2 & \hdots & 2 \\
\vdots  & \vdots &  & \vdots \\
2^{2k-1} & 2^{2k-1} & \hdots & 2^{2k-1} 
\end{bmatrix}
-
\begin{bmatrix}
1 & 0 & \hdots & 0 & 0 \\
2 & 1 & \ddots &  & 0 \\
\vdots  & \vdots & \ddots & \ddots & \vdots \\
2^{2k-2} & 2^{2k-3}& \hdots & 1& 0 \\
2^{2k-1} & 2^{2k-2} & \hdots & 2 & 1
\end{bmatrix}.
$$
\begin{proof}
Using the identity $t_{n+2k+1} = 2t_{n+2k} - t_{n}$, one can show by induction that for any $0 \leq i \leq 2k-1$
\begin{align} \label{eq:matrix_rows}
t_{n +2k+1+i} = 2^{i+1} t_{n+ 2k} - \sum_{j=0}^{i} 2^{i-j} t_{n+j} = 2 \cdot 2^i \sum_{j=0}^{2k-1} t_{n+j}- \sum_{j=0}^{i} 2^{i-j} t_{n+j}.
\end{align}
We know that $D = C^{2k+1}$ is the only matrix satisfying $T_{n+2k+1} = DT_n$ for all $n \geq 0$. Thus, for each $0 \leq i \leq 2k-1$ the coefficients at $t_{n+j}$ in (\ref{eq:matrix_rows}) correspond to the $i$-th row of $C^{2k+1}$ (counting from 0).
\end{proof}
\end{lem}

We are also going to need two standard identities involving binomial coefficients. For the convenience of the reader we include the proof.

\begin{lem}\label{binomial_formulas}
For all positive integers $m, w$ we have

{\rm (a)} $\sum\limits_{i=0}^{w} \binom{m+i}{m} = \binom{m+w+1}{m+1}$,

{\rm (b)} $\sum\limits_{i=0}^{w} \binom{m+i}{m} 2^i = (-1)^{m+1} + 2^{w+1} \sum\limits_{j=0}^{m} \binom{m+w+1}{m-j} (-2)^j.$

\begin{proof}
For any fixed $m \geq 0$ the formula (a) follows easily from induction on $w$.

To prove (b) take any $w,m \geq 0$ and consider the function
$$
f(x) = \frac{1}{m!}\sum\limits_{i=0}^{w} x^{i+m} = \frac{1}{m!} \; \frac{x^{m+w+1}-x^m}{x-1} 
$$
for $x \neq 1$.
It is easy to see that the left side of (b) is equal to $f^{(m)}(x)$ evaluated at $x=2$. Applying the Leibniz formula we get
\begin{align*}
f^{(m)}(x) &= \frac{1}{m!} \sum_{j=0}^{m} \binom{m}{j}  \left(\frac{1}{x-1}\right)^{(j)}  \left(x^{m+w+1}-x^m \right)^{(m-j)}  \\
&=\sum_{j=0}^{m} \frac{(-x)^{j}}{(x-1)^{j+1}}  \left [  \binom{m+w+1}{m-j} x^{w+1} - a\binom{m}{j}  \right ]  \\
&= x^{w+1} \sum_{j=0}^{m} \binom{m+w+1}{n-j} \frac{(-x)^{j}}{(x-1)^{j+1}} +\left( \frac{-1}{x-1} \right)^{m+1}.
\end{align*}
Substituting $x=2$ we get the desired result.
\end{proof}
\end{lem}

The following lemma gives us an easily computable expression for $2k$ subsequent terms $t_n$.
\begin{lem} \label{lem:easy_expression}
Define the following column vectors in $\mathbb{N}^{2k}$:
\begin{align*}
w =  \begin{bmatrix}
1 & 1 & \hdots & 1
\end{bmatrix}^T,
\quad
v_m = \begin{bmatrix}
\binom{m}{m} & \binom{m+1}{m} \cdot 2^1 & \binom{m+2}{m} \cdot 2^2 &
\hdots & \binom{m+2k-1}{m} 2^{2k-1}
\end{bmatrix}^T
\end{align*}
for $m \geq 0$. Then for any $m \geq 1$ we have
\begin{align} \label{eq:whole_push}
T_{m(2k+1)} \equiv w + (-1)^{m+1} \cdot 4(k-1) \sum_{i=0}^{m-1}v_i + (-1)^{m+1} v_{m-1} \pmod{2^{2k + 1}}.
\end{align}

\begin{proof}
Using the form of $C^{2k+1}$ given in Lemma \ref{lem:pushing_the_sequence}, it is easy to see that
\begin{align} \label{eq:first_push}
C^{2k+1} w = 2 \cdot 2k v_0 - (2v_0 - w) = w + 2(2k-1)v_0.
\end{align}
Now fix $m \geq 0$.
Applying Lemma \ref{binomial_formulas} to each coordinate gives us

\begin{align} \label{eq:later_push}
C^{2k+1} v_m \equiv 2\cdot(-1)^{m+1} v_0 - v_{m+1} \pmod{2^{2k+1}}.
\end{align}
One can check that
$$
C^{2k+1} T_0 = w + 2(2k-1)v_0 - v_0 = w + 4(k-1)v_0 + v_0,
$$
so (\ref{eq:whole_push}) is true for $m=1$. Now assume that (\ref{eq:whole_push}) holds for some $m \geq 1$. Using (\ref{eq:first_push}) and (\ref{eq:later_push}) we get
\begin{align*}
T_{(m+1)(2k+1)} &= C^{2k+1} T_{m(2k+1)} \\
&\equiv C^{2k+1} \left[w + (-1)^{m+1} \cdot 4(k-1) \sum_{i=0}^{m-1}v_i + (-1)^{m+1} v_{m-1}\right] \\
&\equiv w + 2(2k-1)v_0 + (-1)^{m+1} \cdot 4(k-1) \sum_{i=0}^{m-1} \left[2 \cdot (-1)^{i+1} v_0 - v_{i+1}\right] \\
&\qquad +(-1)^{m+1} \left[ 2 \cdot (-1)^{m} v_0 - v_{m}\right]  \\
&\equiv w + \left[4(k-1) - 8(k-1) \epsilon_m\right]v_0 + (-1)^{m} \cdot 4(k-1) \sum_{i=1}^{m} v_{i} \\
&\qquad + (-1)^{m} v_m \pmod{2^{2k+1}},
\end{align*}
where $\epsilon_m$ is equal to $m$ modulo 2. Thus, the coefficient at $v_0$ is equal to $(-1)^m \cdot 4(k-1)$, so we can incorporate it into the sum. Finally, we obtain
$$
T_{(m+1)(2k+1)}
\equiv 
w + (-1)^{m} \cdot 4(k-1) \sum_{i=0}^{m} v_{i} + (-1)^{m} v_m \pmod{2^{2k+1}}.
$$
\end{proof}
\end{lem}

We are now ready to prove Theorem \ref{thm:2-adic_order}.

\begin{proof}[Proof (of Theorem \ref{thm:2-adic_order})]
The term $t_n$ is even iff $n$ is divisible by $2k+1$, which proves that $\nu_2(t_n) = 0$ for $n \equiv 1,2,...,2k \pmod{2k+1}$. 
Observe that if $k \geq 2$, then $2 \nu_2(k-1) + 5 \leq 2k+1$, so by Lemma \ref{lem:pushing_the_sequence} for $m \geq 1$ we have
$$
T_{m(2k+1)} \equiv
w + (-1)^{m+1} \cdot 4(k-1) \sum_{i=0}^{m-1}v_i + (-1)^{m+1} v_{m-1} \pmod{2^{2 \nu_2(k-1) + 5}}.
$$
Looking at the first entry of this vector, we obtain
\begin{align} \label{eq:super_formula}
t_{m(2k+1)} \equiv 1 + (-1)^{m+1} \cdot 4m(k-1) + (-1)^{m+1} \pmod{2^{2 \nu_2(k-1) + 5}}.
\end{align}
Thus, for odd $m$ we get $t_{m(2k+1)} \equiv 2 \pmod{4}$, hence $\nu_2(t_n) = 1$ for $n \equiv 2k+1	\pmod{2(2k+1)}$.

Now let $n = s(2k+1) 2^l$ for odd $s$ and $l \geq 1$, so that $n \equiv 0 \pmod{2(2k+1)}$. We will further split the third case into two subcases, depending whether $l \leq l_0 = \nu_2(k-1) +2$ or $l > l_0$. 
If $l \leq l_0$ then from (\ref{eq:super_formula}) we obtain 
\begin{align} \label{eq:subcase_1}
t_n \equiv s2^{l+2} (k-1) \pmod{2^{2 \nu_2(k-1) + 5}},
\end{align} 
so $\nu_2(t_n) = l+ 2 + \nu_2(k-1) = \nu_2(n) + \nu_2(k-1) + 2$.

We cannot extend the same argument to $l > l_0$ because we only know the congruence modulo $2^{2 \nu_2(k-1) + 5}$. However, substituting $s=1$ and $l = l_0$ in (\ref{eq:subcase_1})  gives us a possible value $a_{l_0,0} = 2(k-1)$, as defined in Lemma \ref{lem:later_terms}. Using Lemma \ref{lem:later_terms} for any $l > l_0$, we get in the first coordinate
$$
t_n \equiv 2^{l+1} a_{l_0,0} \equiv 2^{l+2} (k-1) \pmod{2^{l + \nu_2(k-1) + 3}},
$$
which again yields $\nu_2(t_n) = \nu_2(n) + \nu_2(k-1) + 2$.
\end{proof}


\begin{thebibliography}{1}
\bibitem{Allouche_Shallit} J.-P. Allouche and J. Shallit, The ring of $k$-regular sequences, \emph{Theor. Comput. Sci.} {\bf 98}: 163--197 (1992).
\bibitem{Allouche_Shallit_2} J.-P. Allouche and J. Shallit, The ring of $k$-regular sequences II, \emph{Theor. Comput. Sci.} {\bf 307}: 3--29 (2003).
\bibitem {Lengyel_order}  T. Lengyel, The order of Fibonacci and Lucas numbers, \emph{Fibonacci Quart.} {\bf 33}: 234--239 (1995).  
\bibitem{Lengyel_multi} T. Lengyel, Divisibility properties by multisection, \emph{Fibonacci Quart.} {\bf 41}: 72--79 (2003).
\bibitem{Lengyel_Marques} T. Lengyel and D. Marques, The 2-adic order of the Tribonacci Numbers and the equation $T_n = m!$, \emph{J. Integer Seq.} {\bf 17}: Article 14.10.1 (2014).
\bibitem{Marques} D. Marques, The order of appearance of product of consecutive Fibonacci numbers, \emph{Fibonacci Quart.} {\bf 50}: 132--139 (2012).
\bibitem{Shu_Yao} Z. Shu and J.-Y. Yao, Analytic functions over $\mathbb{Z}_p$ and $p$-regular sequences, \emph{C. R. Math.} {\bf 349}: 947--952 (2011).
\bibitem{Wolfram} D. A. Wolfram, Solving generalized Fibonacci recurrences, \emph{Fibonacci Quart.} {\bf 36}: 129--145 (1998)

\end{thebibliography}
\end{document}